\documentclass[12pt]{article}
\usepackage{a4wide}
\usepackage{theorem}
\usepackage{amssymb}
\usepackage{amsmath}
\usepackage{amsfonts}

\newtheorem{theo}{Theorem}[section]

\newtheorem{prop}[theo]{Proposition}

\newtheorem{lem}[theo]{Lemma}
\newtheorem{rem}[theo]{Remark}

\newtheorem{coro}[theo]{Corollary}
\newtheorem{defi}[theo]{Definition}
\newtheorem{exam}[theo]{Example}

\newcommand{\biindice}[3]%
{\renewcommand{\arraystretch}{0.5}
\begin{array}[t]{c}
#1\\
{\scriptscriptstyle #2}\\
{\scriptscriptstyle #3}
\end{array}
\renewcommand{\arraystretch}{1}}

\newcommand{\geqs}[1]{\arraycolsep0.1pt \renewcommand{\arraystretch}{0.5}
\begin{array}[t]{c}
\geq\\
#1
\end{array}
\arraycolsep5pt \renewcommand{\arraystretch}{1}
}

\newcommand{\leqs}[1]{\arraycolsep0.1pt \renewcommand{\arraystretch}{0.5}
\begin{array}[t]{c}
\leq\\
#1
\end{array}
\arraycolsep5pt \renewcommand{\arraystretch}{1}
}

\newcommand{\dR}{I \! \! R}

\title{ON THE FU\v CIK SPECTRUM WITH INDEFINITE WEIGHTS}
\author{~~\\[10mm]M. ALIF and J.-P. GOSSEZ\\
D\'ep. Math\'ematique, C.P. 214,\\
Universit\'e Libre de Bruxelles\\
1050 Bruxelles, Belgium\\
e-mail:\\
gossez@ulb.ac.be\\
malif@ictp.trieste.it
}
\date{}

\begin{document}

\maketitle

\newpage

\section{Introduction}

This paper is partly concerned with the one-dimensional asymmetric 
problem with weight
\begin{equation}
    \left\{
    \begin{array}{l}
	Lu = m(t) (au^{+}-bu^{-}) \mbox{ in } ]T_{1},T_{2}[,\\
	u(T_{1}) = u(T_{2})=0.
    \end{array}
\right.
    \label{eq:1.1}
\end{equation}
Here $Lu:=-(p(t)u')'+q(t)u$, $p,q$ and $m \in C[T_{1},T_{2}]$, $p(t) > 
0$ on $[T_{1},T_{2}]$, $q(t) \geq 0$ on $[T_{1},T_{2}]$, $m(t) \not 
\equiv 0$ and  $u^{\pm}:=\max\{\pm u,0\}$. The associated Fu\v cik 
spectrum is defined as the set $\Sigma$ of those $(a,b) \in \dR^{2}$ 
such that (\ref{eq:1.1}) has a nontrivial solution $u$.

\vspace{5mm}

The description of this spectrum $\Sigma$ is classical and explicit when $Lu=-u''$ 
and there is no weight, i.e. $m(t) \equiv 1$ (cf. \cite{Da}, 
\cite{Fu}). The same general picture for $\Sigma$ remains valid when 
$L$ is as above and $m(t) > 0$ on $]T_{1}, T_{2}[$ (cf. \cite{Dr}, 
\cite{Ca}, \cite{Ry})~: $\Sigma$ is made of the two lines $\dR 
\times \lambda_{1}^{m}$ and $\lambda_{1}^{m} \times \dR$ together with 
a sequence of hyperbolic like curves in $\dR^{+}\times \dR^{+}$ 
passing through $(\lambda_{k}^{m}, \lambda_{k}^{m}$), $k \geq 2$; one 
or two such curves emanate from each $(\lambda_{k}^{m}, 
\lambda_{k}^{m})$, and the corresponding solutions of (\ref{eq:1.1}) 
along these curves have exactly $k-1$ zeros in $]T_{1},T_{2}[$. 
Here $(0 <) \lambda_{1}^{m} < \lambda_{2}^{m} < \ldots \to + \infty$ 
denotes the sequence of eigenvalues of the associated linear problem
\begin{equation}
    \left\{
    \begin{array}{l}
    Lu  = \lambda m(t) u \mbox{ in } ]T_{1},T_{2}[,\\
    u(T_{1}) = u(T_{2}) = 0.
\end{array}
\right.
    \label{eq:1.2}
\end{equation}

One of our purposes in this paper is to investigate the situation where 
the weight function $m(t)$ in (\ref{eq:1.1}) {\sl changes sign} in 
$]T_{1},T_{2}[$. In that case it is well-known that the eigenvalues 
in (\ref{eq:1.2}) form a double sequence~:
$- \infty \leftarrow \ldots < \lambda_{-2}^{m} < \lambda_{-1}^{m} < 
(0) < \lambda_{1}^{m} < \lambda_{2}^{m} < \ldots \to + \infty.
$
A natural conjecture is then that $\Sigma$ should be made of the four 
trivial lines $\dR \times \lambda_{1}^{m}$, $\lambda_{1}^{m} \times 
\dR$, $\dR \times \lambda_{-1}^{m}$, $\lambda_{-1}^{m} \times \dR$ 
together with a sequence of hyperbolic like curves in $\dR^{+} \times 
\dR^{+}$ passing through $(\lambda_k	^m	, \lambda_k	^m	)$, $k \geq 2$, and a sequence 
of hyperbolic like curves in $\dR^{-} \times \dR^{-}$ passing through
$(\lambda_{-k}^{m},\lambda_{-k}^{m})$,
$k \geq 2$. We indeed find such a double sequence of curves but, rather 
surprisingly, other hyperbolic like curves {\sl always} appear in the 
other two quadrants $\dR^{+} \times \dR^{-}$ and $\dR^{-} \times 
\dR^{+}$. The number of these additional curves depends on 
the number of changes of sign of the weight $m(t)$: if $m(t)$ changes sign $N(=1,2,\ldots,+\infty)$ times
in
$]T_{1},  T_{2}[$ (cf. definition 4.3), then one gets exactly $(2N-1)$ curves in 
$\dR^{+} \times \dR^{-}$ and also in $\dR^{-} \times \dR^{+}$. (Note 
that some of these curves may be double and are then counted for two, 
cf. remark 4.5). These additional curves as well as those passing 
through $(\lambda_{k}^{m}, \lambda_{k}^{m})$ and $(\lambda_{-k}	^m	, \lambda_{-k}	^m	)$ can be classified 
according to the number of zeros of the corresponding solutions of 
(\ref{eq:1.1}).

\vspace{3mm}

More generally we consider the problem
\begin{equation}
    \left\{
    \begin{array}{l}
	Lu = am(t)u^{+}-bn(t)u^{-} \mbox{ in } ]T_{1},T_{2}[,\\
	u(T_{1}) = u(T_{2}) = 0
    \end{array}
\right.
    \label{eq:1.3}
\end{equation}
which involves {\sl two} weights $m,n \in C[T_{1},T_{2}]$ with $m(t)$ 
and $n(t) \not \equiv 0$. Let again $\Sigma$ denote the corresponding 
Fu\v cik spectrum. Note that in this case, making $a=b$ in 
(\ref{eq:1.3}) does not lead to any obvious elements in $\Sigma$.  
Assuming for instance that both $m(t)$ and $n(t)$ change sign in 
$]T_{1}, T_{2}[$, we show that beside a trivial part consisting of the 
4 lines $\lambda_{1}^{m} \times \dR$, $\lambda_{-1}^{m} \times \dR$, 
$\dR \times \lambda_{1}^{n}$ and $\dR \times \lambda_{-1}^{n}$, 
$\Sigma$ is made in {\sl each} quadrant of a (non zero) odd or 
infinite number of hyperbolic like curves. These curves can again be 
classified according to the number of zeros of the corresponding 
solutions of (\ref{eq:1.3}). We also show that all cases can 
effectively happen with respect to the numbers of these curves : 
given $K, L, M, N \in \{0,1,2,\ldots,+\infty\}$, there exist weights 
$m(t)$ and $n(t)$ such that $\Sigma$ exactly contains $(2K+1), 
(2L+1), (2M+1)$ and $(2N+1)$ hyperbolic like curves in $\dR^{+} 
\times \dR^{+}$, $\dR^{-} \times \dR^{-}$, $\dR^{+} \times \dR^{-}$ 
and $\dR^{-} \times \dR^{+}$ respectively. (Note that here as before 
some of these curves may be double and are then counted for two, cf. 
remark 3.9).

\vspace{5mm}

Section 3 deals with the two weights problem (\ref{eq:1.3}) and 
section 4 with the one weight problem (\ref{eq:1.1}). In section 5 we 
investigate for (\ref{eq:1.3}) the asymptotic behaviour of the first 
hyperbolic like curves of $\Sigma$ (i.e. those which lie the closest 
to the trivial horizontal and vertical lines). It is known that if 
$m(t)$ and $n(t)$ are $>0$ in $]T_{1},T_{2}[$, then these first 
curves are asymptotic to the line $\dR \times \lambda_{1}^{n}$ as $a 
\to + \infty$ and to the line $\lambda_{1}^{m} \times \dR$ as $b 
\to + \infty$ (cf. \cite{Ry}; cf. also \cite{Dr}, \cite{Ca} in the one 
weight problem). As observed in \cite{DF-Go}, such an asymptotic 
behaviour is closely connected with the nonuniformity of the 
antimaximum principle. It turns out that this asymptotic behaviour may 
be affected by the presence of more general weights. We show in 
particular that if $m(t)$ and $n(t)$ have compact support in $]T_{1}, 
T_{2}[$, then none of the first curves is asymptotic on any side to 
the trivial horizontal and vertical lines; the converse implication 
is also true. In remark \ref{rem:5.6}, we briefly comment on the meaning of this result in the context of 
the antimaximum principle.

\vspace{3mm}

Our approach is based on the shooting method and in section 2, which 
has a preliminary character, we investigate various properties of the 
``zero-function''. Given a nontrivial solution $u$ of the linear 
equation $Lu=am(t)u$, this function sends one zero of $u$ onto the 
following zero of $u$.

\vspace{5mm}

We finally mention that results analogous to those in the present 
paper can also be established for the Neumann problem (cf. \cite{Al1}, \cite{Al2}). 
Moreover the study of the Fu\v cik spectrum with indefinite weights in the P.D.E. case has been initiated
recently in \cite{star}.

\section{The Zero-function}
\setcounter{equation}{0}

In this section we consider the zero-function for the linear 
equation $Lu=am(t)u$. Here $L$ is as in the introduction and $m \in 
C[T_{1}, T_{2}]$, $m(t) \not \equiv 0$. It will be convenient to 
extend the coefficients of $L$ and the weight $m$ from $[T_{1}, 
T_{2}]$ to the whole of $\dR$, preserving continuity as well as the 
inequalities $p_{1} \leq p(t) \leq p_{2}$, $0 \leq q(t) \leq q_{2}$ 
and $m_{1} \leq m(t) \leq m_{2}$ for some constants $p_{1}, 
p_{2},q_{2},m_{1},m_{2}$ with $p_{1} > 0$. We will assume below that 
such an extension has been carried out.

\vspace{5mm}

Writing the equation $Lu=am(t)u$ as a first order system, it follows 
from standard ODE results (cf. e.g. ch. 5 of \cite{Ha}) that for each $s 
\in \dR$, the initial value problem
\begin{equation}
    \left\{
    \begin{array}{l}
	Lu=am(t)u \mbox{ in } \dR,\\
	u(s) = 0, \quad u'(s) = 1/p(s)
    \end{array}
\right.
    \label{eq:2.1}
\end{equation}
has a unique solution $u=u(t)=u(t;a,s)$. This solution is a $C^{1}$ 
function of $(t,a,s) \in \dR \times \dR \times \dR$; moreover $u$ 
and $(pu')$ have second mixed derivatives with respect to $t,a$, and 
also with respect to $t,s$, which commute and are continuous.

\begin{defi}\label{defi:2.1}
    {\rm The zero-function $\varphi_{}$ is defined by
    $$
    \varphi_{a}(s) : = \min \{ t \in \dR : t > s \mbox{ and } 
    u(t;a,s)=0\},
    $$
    with $\varphi_a(s)=+\infty$ in case $u(t)$ does not vanish for any 
    $t > s$.
}
\end{defi}

Thus $\varphi_{a}$ sends  $s$ onto the following zero of the solution $u$ 
of (\ref{eq:2.1}). Since the zeros of $u$ are isolated, this 
definition makes sense. Clearly $\varphi_{a}(s) > s$. Note that $\varphi_a	(s)$ is also the first zero following $s$ of any
nontrivial solution of $Lu =am(t)u$, $u(s)=0$.

\vspace{5mm}

Some monotonicity and regularity properties of this function are 
collected in the following two lemmas. We call $D$ the domain of 
$\varphi$, i.e. $D:=\{(a,s) \in \dR^{2} : \varphi_{a}(s) < + \infty\}$.

\begin{lem}
    \label{lem:2.2}
    (i) For each $a, \varphi_{a}(s)$ is increasing with respect to 
    $s$, strictly on $D$. 
    (ii) For each $s$, $\varphi_{a}(s)$ is 
    decreasing (resp. increasing) with respect to $a$ for $a \geq 0$ 
    (resp. for $a \leq 0$), strictly on $D$.
\end{lem}

{\bf Proof.} These monotonicity properties easily follow from the Sturm 
comparison theorem as given e.g. in ch. 11 of \cite{Ha}. To apply 
this theorem in part (ii), it is useful to write first the equation 
$Lu=am(t)u$ as $-((p(t)/a)u')'+(q(t)/a)u=m(t)u$. Q.E.D.

\begin{lem}\label{lem:2.3}
    (i) $D$ is open and $\varphi : D \to \dR$ is of class $C^{1}$.
    (ii) $\partial \varphi_{a}(s)/\partial s  > 0$ for $(a,s) \in D$. 
    (iii) $\partial \varphi_a(s)/\partial a < 0$ for $(a,s)  \in D$ 
    with $a > 0$. (iv) $\partial \varphi_a(s)/\partial a > 0$ for 
    $(a,s) \in D$ with $a < 0$.
\end{lem}

{\bf Proof.} We start with (i). Let $(a_{0},s_{0}) \in D$. Applying 
the implicit function theorem, we get open neighbourhoods $U$ of 
$(a_{0},s_{0})$ and $V$ of $\varphi_{a_{0}}(s_{0})$ and a $C^{1}$ 
function $\tilde{\varphi}~: U \to V$ such that 
$$
\tilde{\varphi}(a,s) = t \mbox{ iff } (t,a,s) \in V \times U \mbox{ 
and } u(t; a,s)=0.
$$
Since $\tilde{\varphi}(a_{0},s_{0})=\varphi_{a_{0}}(s_{0})>s_{0}$, 
reducing $U$ if necessary, we can assume $\tilde{\varphi}(a,s) > s$ 
on $U$.
The definition of $\varphi_{a}(s)$ then implies $\varphi_{a}(s) \leq \tilde{\varphi}(a,s)$. 
Consequently $U \subset D$ and $D$ is open. We will now show that 
$\varphi \equiv \tilde{\varphi}$ near $(a_{0},s_{0})$, which will 
conclude the proof of (i). Assume by contradiction the existence of a 
sequence $(a_{i},s_{i}) \to (a_{0},s_{0})$ such that 
\begin{equation}
    s_{i} < \varphi_{a_{i}}(s_{i}) < \tilde{\varphi}(a_{i},s_{i}).
    \label{eq:2.2}
\end{equation}
For a subsequence, $\varphi_{a_{i}}(s_{i}) \to $ some $r$ which 
satisfies $r \in [s_{0}, \varphi_{a_{0}}(s_{0})]$ and $u(r; a_{0}, 
s_{0})=0$. This implies $r=s_{0}$ or $r=\varphi_{a_{0}}(s_{0})$. By 
(\ref{eq:2.2}) we can apply Rolle's theorem to $u(t; a_{i},s_{i})$ 
between either $s_{i}$ and $\varphi_{a_{i}}(s_{i})$, or 
$\varphi_{a_{i}}(s_{i})$ and $\tilde{\varphi}(a_{i},s_{i})$. We deduce 
that either $u'(s_{0}; a_{0},s_{0})=0$ or $u'(\varphi_{a_{0}}(s_{0}); 
a_{0}, s_{0})=0$, a contradiction.

We now turn to the proof of (ii), (iii), (iv). These properties easily follow by 
derivating the relation $u(\varphi_{a}(s);a,s) \equiv 0$ once we know 
that 
\begin{equation}
    \left[\partial u(t;a,s)/\partial t\right]_{t=\varphi_{a}(s)} < 0,
    \label{eq:2.3}
\end{equation}
\begin{equation}
     \left[\partial u(t;a,s)/\partial s\right]_{t=\varphi_{a}(s)} > 0,
    \label{eq:2.4}
\end{equation}
\begin{equation}
     \left[\partial u(t;a,s)/\partial a\right]_{t=\varphi_{a}(s)} < 0 \mbox{ if } a > 0 \mbox{ 
    and } > 0 \mbox{ if } a < 0.
    \label{eq:2.5}
\end{equation}
(\ref{eq:2.3}) is clear from the definition of $\varphi_{a}(s)$. To 
prove (\ref{eq:2.4}) we define $v(t) := \partial u(t; a,s)/\partial s$.
Derivating the relation $u(s;a,s) \equiv 0$ with respect to $s$, we 
get $v(s) = -1/p(s) \neq 0$.
Derivating the equation (\ref{eq:2.1}) with respect to $s$, we get 
$Lv=am(t)v$. So $v(\varphi_{a}(s)) \neq 0$ and, since $v(s) < 0$, the Sturm 
comparison theorem implies $v(\varphi_{a}(s)) > 0$, i.e. (\ref{eq:2.4}). 
To prove (\ref{eq:2.5}) we define $w(t) := \partial u(t;a;s)/\partial a)$. Since $u(s;a,s) \equiv 0$, we
get $w(s)=0$. Derivating the  equation in (\ref{eq:2.1}) with respect to $a$, we get 
$Lw=am(t)w+(Lu)/a$. Multiplying by $u$ and integrating from $s$ to 
$\varphi_{a}(s)$, we obtain, after two integrations by part,
$$
p(\varphi_{a}(s))w(\varphi_{a}(s))u'(\varphi_{a}(s))=\frac{1}{a} 
\int_{a}^{\varphi_{a}(s)}p(t)(u')^{2}+\frac{1}{a}\int_{a}^{\varphi_{a}(s)}q(t)u^{2},
$$
which implies (\ref{eq:2.5}). Q.E.D.

\vspace{5mm}

The behaviour of $\varphi_{a}(s)$ as $a \to \pm \infty$ will be 
important in our study. Let us define
\begin{equation}
    \begin{array}{l}
	\alpha_{s}^{>} := \inf \{t > s : m(t) > 0\},\\
	\alpha_{s}^{<} : = \inf \{t > s : m(t) < 0\},
    \end{array}
    \label{eq:2.6}
\end{equation}
with the usual convention that $\inf \phi = + \infty$. In  the 
simplest cases, $\alpha_{s}^{>}$ (resp. $\alpha_{s}^{<}$) represents 
the lower bound of the first positive (resp. negative) bump of $m(t)$ 
situated at the right of $s$. Clearly $\alpha_{s}^{>} \geq s$ and 
$\alpha_{s}^{<} \geq s$.

\begin{lem}\label{lem:2.4}
    For each $s \in \dR$, $\varphi_{a}(s) \to \alpha_{s}^{>}$ as $a 
    \to +\infty$ and $\varphi_a	(s) \to \alpha_s	^<	$ as $a \to -\infty$.
\end{lem}

{\bf Proof.} We will consider only the case $a \to + \infty$ (the 
other case can be reduced to this one by considering $-m$). We will 
also assume $\alpha_{s}^{>} < +\infty$ (otherwise the conclusion is 
trivial). 

Assume first that $s$ is such that $m(s) > 0$. Given $\epsilon > 0$ 
sufficiently small, $m$ is $> 0$ on $[s,s+\epsilon]$, and 
consequently there exists $a_{\epsilon}$ such that $am(t) \geq 
p_{2}(\pi/\epsilon)^{2}+q_{2}$ on $[s,s+\epsilon]$ for all $a \geq 
a_{\epsilon}$, where $p_{2}$ and $q_{2}$ are defined at the beginning 
of section 2. We will compare on $[s,s+\epsilon]$ our equation 
$Lu=am(t)u$ with the equation
$$
-p_{2}v''=p_{2}(\pi/\epsilon)^{2}v.
$$
Since $\sin(\pi(t-s)/\epsilon)$ is a solution of the latter which 
vanishes at $s$ and $s + \epsilon$, the Sturm comparison theorem 
implies that $s < \varphi_{a}(s) \leq s + \epsilon$. The conclusion of 
the lemma then follows since $\alpha_{s}^{>} = s$ in the case under consideration.

We now turn to the general case. Let $\epsilon > 0$. By the definition of $\alpha_{s}^{>}$, we can 
find $s_{\epsilon} \in ]\alpha_{s}^{>},\alpha_{s}^{>}+\epsilon[$ such 
that $m(s_{\epsilon}) > 0$. The case already treated then implies that for $a$ sufficiently large,
\begin{equation}
    \varphi_{a}(s) < \varphi_{a}(s_{\epsilon}) \leq s_{\epsilon} + 
    \epsilon \leq \alpha_{s}^{<} + 2\epsilon.
    \label{eq:2.7}
\end{equation}
On the other hand one has
\begin{equation}
    \alpha_{s}^{>} < \varphi_{a}(s).
    \label{eq:2.8}
\end{equation}
Indeed, if (\ref{eq:2.8}) does not hold, then $m$ is $\leq 0$ on 
$[s,\varphi_{a}(s)]$. Comparing on this interval our equation $Lu=amu$ 
with the equation $-p_{1} v''=0v$, where $p_{1}$ is defined at the 
beginning of section 2, we deduce from the Sturm comparison theorem 
that any solution $v$ of the latter equation must have a zero in 
$[s,\varphi_{a}(s)]$, which is clearly false. Combining 
(\ref{eq:2.7}) and (\ref{eq:2.8}) finally yields the conclusion of the 
lemma. Q. E. D.

\begin{lem}\label{lem:2.5}
     For each $s \in \dR$, $\varphi_{a}(s) \to + \infty$ as $a \to 0$.
\end{lem}

{\bf Proof.} We consider only the case $a \geq 0$ (the case $a \leq 0$ 
can be reduced to this one by considering $-m$). Let $R > 0$.
There exists $a_{R} > 0$ such that $am(t) \leq p_{1}(\pi/R)^{2}$ for 
$0 \leq a \leq a_{R}$, where $p_{1}$ is defined at the beginning of 
section 2. We compare on $[s,s+R]$ our equation $Lu=am(t)u$ with the 
equation
$$
-p_{1}v^{''}=p_{1}(\pi/R)^{2}v.
$$
Since $\sin(\pi(t-s)/R)$ is a solution of the latter which has $s$ 
and $s+R$ as consecutive zeros, the Sturm comparison theorem 
implies $\varphi_{a}(s) \geq s+R$. Q. E. D.

\vspace{3mm}

To conclude this section we return to our given interval 
$[T_{1},T_{2}]$. Clearly the restriction of the function 
$\varphi_{a}$ to the set $\{s \in [T_{1},T_{2}]; \varphi_{a}(s) \in 
[T_{1},T_{2}]\}$ does not depend on the extension of the datas carried 
out at the beginning of this section.

\vspace{3mm}

If $m^{+}(t) \not \equiv 0$ on $]T_{1}, T_{2}[$, then 
$\alpha_{T_{1}}^{>}$ defined in (\ref{eq:2.6}) is $< T_{2}$, and it 
follows from lemmas \ref{lem:2.2}, \ref{lem:2.3}, \ref{lem:2.4} and 
\ref{lem:2.5} that $\varphi_{a}(T_{1})=T_{2}$ for exactly one value 
of $a > 0$. Clearly this value of $a$ is equal to the first positive 
eigenvalue $\lambda_{1}^{m}$ of (\ref{eq:1.2}). Similarly, if 
$m^{-}(t) \not \equiv 0$ on $]T_{1},T_{2}[$, then 
$\varphi_{a}(T_{1})=T_{2}$ for exactly one value of $a < 0$, which is 
equal to $\lambda_{-1}^{m}$.

\vspace{3mm}

The zero-function $\varphi_{a}$ clearly depends on the weight $m(t)$ 
and we will from now on denote it by $\varphi_{a}^{m}$.

\section{The two weights problem}
\setcounter{equation}{0}

In this section we consider problem (\ref{eq:1.3}), with $L$ as in 
the introduction  and $m,n \in C[T_{1}, T_{2}]$, $m(t)$ and $n(t) 
\not \equiv 0$.

\vspace{5mm}

Since the zeros of a nontrivial solution of (\ref{eq:1.3}) are 
isolated, one can classify these solutions 
according to the number of their zeros. This yields the following 
description of the Fu\v cik spectrum $\Sigma$~:

\begin{prop}\label{prop:3.1}
   We have 
   $$
   \Sigma = C_{1}^{>} \cup C_{1}^{<} \cup C_{2}^{>} \cup C_{2}^{<} 
   \cup C_{3}^{>} \cup C_{3}^{<} \cup \ldots
   $$
   where
   \begin{center}
      $C_{k}^{>}$ (resp. $C_{k}^{<})=\{(a,b) \in \dR^{2}$; 
      (\ref{eq:1.3})  has a solution\\
      $u$ with $k-1$ zeros in $]T_{1},T_{2}[$ and $u'(T_{1}) > 0$ 
      (resp. $u' (T_{1}) < 0\}$.
   \end{center}
   Moreover, if $(a,b) \in C_{k}^{>}$ (resp. $C_{k}^{<}$), then all 
   nontrivial solutions of (\ref{eq:1.3}) with $u'(T_{1}) > 0$ 
   (resp. $u'(T_{1}) < 0$) are multiple one of the other.
\end{prop}

These sets $C_{k}^{>}$ and $C_{k}^{<}$ can themselves be described in 
terms of the zero-functions $\varphi_{a}^{m}$ and $\varphi_{b}^{n}$~:
\begin{eqnarray}
    C_{1}^{>} & = & \{(a,b) \in \dR^{2} : 
    \varphi_{a}^{m}(T_{1})=T_{2}\},
    \nonumber  \\
    C_{2}^{>} & = & \{(a,b) \in \dR^{2} : 
    \varphi_{b}^{n}(\varphi_{a}^{m}(T_{1}))=T_{2}\},
    \nonumber  \\
    C_{3}^{>} & = & \{(a,b) \in \dR^{2} : 
    \varphi_{a}^{m}(\varphi_{b}^{n}(\varphi_{a}^{m}(T_{1})))=T_{2}\},\ldots
    \nonumber\\
    ~~~ &\,& \label{eq:3.1}\\
    C_{1}^{<} & = & \{(a,b) \in \dR^{2} : 
    \varphi_{b}^{n}(T_{1})=T_{2}\},
    \nonumber  \\
    C_{2}^{<} & = & \{(a,b) \in \dR^{2} : 
    \varphi_{a}^{m}(\varphi_{b}^{n}(T_{1}))=T_{2}\},
    \nonumber  \\
    C_{3}^{<} & = & \{(a,b) \in \dR^{2} : 
    \varphi_{b}^{n}(\varphi_{a}^{m}(\varphi_{b}^{n}(T_{1})))=T_{2}\},\ldots
    \nonumber
\end{eqnarray}
Of course, as we will see later, some of these sets may be empty. Note 
that $C_{1}^{>}$ is made of one or two vertical lines 
$\lambda_{1}^{m} \times \dR$ and $\lambda_{-1}^{m} \times \dR$ 
(depending on whether $m(t)$ is $\geq 0$, $\leq 0$ or changes sign). 
Similarly for $C_{1}^{<}$ with the horizontal lines $\dR \times 
\lambda_{1}^{n}$ and $\dR \times \lambda_{-1}^{n}$.
It will be convenient to denote by $\Sigma^{*}$ the set $\Sigma$ 
without these (two, three or four) trivial lines.

\vspace{5mm}

The monotonicity properties of the zero-function imply that if for 
instance $m(t)$ and $n(t)$ change sign, then $\Sigma^{*}$ is contained 
in the four quadrants
$$
(]\lambda_1^m,+\infty[\times]\lambda_1^n,+\infty[) \cup (]-\infty, \lambda_{-1}^{m}[ \, \times \, 
]-\infty,\lambda_{-1}^{n}[)
$$ 
$$\cup (]-\infty,\lambda_{-1}^m[ \, \times \, ]\lambda_1^n,+\infty[)
\cup  (]\lambda_{1}^{m}, + \infty [\,  \times \, ]-\infty, \lambda_{-1}^n[).
$$
In the case where for instance $m(t)$ changes sign but $n(t)$ is, say, $\leqs{\not \equiv} 0$, then
$\Sigma^*	$ is contained in the two quadrants
$$
(]-\infty,\lambda_{-1}^m[ \, \times \, ]-\infty,\lambda_{-1}^n[)\cup (]\lambda_1^m,+\infty[\, \times \,
]-\infty,\lambda_{-1}^n[).
$$  
If $m(t)$ and $n(t)$ do not change sign, then $\Sigma^{*}$ is 
contained in one quadrant.

\begin{rem}\label{rem:3.2}
    {\rm The sets $C_{k}^{>}$ and $C_{k}^{<}$ depend on the 
    weights~: $C_{k}^{>}=C_{k}^{>}(m,n)$ and $C_k^<=C_{k}^{<}(m,n)$. Clearly
    $$
    (a,b) \in C_{k}^{>} (m,n) \Leftrightarrow (a,-b) \in 
    C_{k}^{>}(m,-n) \Leftrightarrow \ldots
    $$
    and similarly for $C_{k}^{<}$. It follows that by changing the 
    sign of the weights, the study of the Fu\v cik spectrum in $\dR 
    \times \dR$ for the weights $(m,n)$ can be reduced to the study 
    of the intersection with $\dR^{+} \times \dR^{+}$ of the Fu\v cik 
    spectrum for the weights $(\pm m, \pm n)$. We will denote below 
    $C_{k}^{>} \cap (\dR^{+} \times \dR^{+})$ by $C_{k}^{>++}$, and 
    similarly for $C_{k}^{<}$ and for the other quadrants. So, for 
    instance, $C_{k}^{<-+} := C_{k}^{<} \cap (\dR^{-} \times 
    \dR^{+})$ and one has 
    $$
    C_{k}^{<-+}(m,n) = \{(a,b): (-a,b) \in C_{k}^{<++}(-m,n)\}.
    $$
}
\end{rem}

As a consequence of this remark, we will often limit ourselves 
below to the study of the intersection of $\Sigma^{*}$ with 
$\dR^{+} \times \dR^{+}$.
    
\vspace{5mm}
    
It will be convenient to introduce a notation for the functions 
which appear in the description (\ref{eq:3.1}) of the sets 
$C_{k}^{>}, C_{k}^{<}$. We will denote by $\Phi_{k}^{>}(a,b)$ 
(resp. $\Phi_{k}^{<}(a,b)$) the composition of $k$ alternating 
functions $\varphi_{a}^{m}$ and $\varphi_{b}^{n}$ starting with 
$\varphi_{a}^{m}$ (resp. $\varphi_{b}^{n}$). So, for instance, 
$$
\Phi_{3}^{>}(a,b)(s)=\varphi_{a}^{m}(\varphi_{b}^{n}(\varphi_{a}^{m}(s))), \, 
\Phi_{3}^{<}(a,b)(s):=\varphi_{b}^{n}(\varphi_{a}^{m}(\varphi_{b}^{n}(s))).
$$
When considering these functions, we will generally assume that, as 
in section 2, the datas have been extended from $[T_{1},T_{2}]$ to the 
whole of $\dR$. So, with this notation
\begin{equation}
    C_{k}^{>++} = \{(a,b) \in \dR^{+} \times \dR^{+} : 
    \Phi_{k}^{>}(a,b)(T_{1})=T_{2}\}
    \label{eq:3.2}
\end{equation}
and similarly for $C_{k}^{<++}$ and for the other quadrants. These 
functions $\Phi_{k}^{>}$ and $\Phi_{k}^{<}$ enjoy properties which 
are rather similar to those of the zero-function.

\vspace{5mm}

The following lemma will be used repeatedly. It characterizes the 
nonemptiness of the sets $C_{k}^{>++}$, $C_{k}^{<++}$ in terms of 
$\Phi_{k}^{>}$, $\Phi_{k}^{<}$ respectively.

\begin{lem}\label{lem:3.3}
    Let $k \geq 2$. Then 
    $$
    C_{k}^{>++} \neq \emptyset \Leftrightarrow \lim_{a,b \to + \infty} 
    \Phi_{k}^{>}(a,b)(T_{1}) < T_{2},
    $$
    $$
    C_{k}^{<++} \neq \emptyset \Leftrightarrow \lim_{a,b \to + \infty} \Phi_{k}^{<} 
    (a,b)(T_{1}) < T_{2}.
    $$
\end{lem}

{\bf Proof.} Note that the limits above exist since, by lemma 2.2, 
$\Phi_{k}^{>}(a,b)(T_{1})$ and $\Phi_{k}^{<}(a,b)(T_{1})$ are 
separately decreasing with respect to $(a,b) \in \dR^{+} \times \dR^{+}$. 
Let us prove the statement relative to $C_{k}^{>++}$ (the other one is 
proved similarly). The implication $\Rightarrow$ follows from 
(\ref{eq:3.2}) and the fact that $\Phi_{k}^{>}(a,b)(T_{1})$ is 
separately strictly decreasing on its domain in $\dR^{+} \times \dR^{+}$. 
The converse implication follows from the continuity of 
$\Phi_{k}^{>}(a,b)(T_{1})$ and the fact that, by lemma 
\ref{lem:2.5}, $\Phi_{k}^{>}(a,b) \to +\infty$ as $a$ or $b \to 0$ Q. 
E. D.

\vspace{3mm}

Our first main result in this section shows that $\Sigma^{*}$ is made 
of hyperbolic like curves of class $C^{1}$.

\begin{theo}\label{theo:3.4}
   Let $k \geq 2$ and assume $C_{k}^{>++}$ nonempty. Then there 
   exist $\nu_{k}^{>} \geq \lambda_{1}^{m}$, $\mu_{k}^{>} \geq 
   \lambda_{1}^{n}$ and a strictly decreasing $C^{1}$ function 
   $f_{k}^{>}$ from $]\nu_{k}^{>}, + \infty[$ onto $]\mu_{k}^{>}, + 
   \infty[$, with $(f_{k}^{>})'(x) < 0$ for $x \in 
   ]\nu_{k}^{>},+\infty[$, such that 
   $$
   C_{k}^{>++}=\{(x,f_{k}^{>}(x)): x \in ]\nu_{k}^{>}, + \infty[\}.
   $$
Similar result for $C_{k}^{<++}$
\end{theo}

{\bf Proof.} Let us prove the result relative to $C_{k}^{>++}$ (the 
other one is proved similarly). Call
$$
I:=\{a > 0 : \exists b > 0 \mbox{ with } 
\Phi_{k}^{>}(a,b)(T_{1})=T_{2}\}.
$$
Since $C_{k}^{>++} \neq \emptyset$, $I \neq \emptyset$. Moreover, by the 
strict monotonicity of $\Phi_{k}^{>}$, for any $a \in I$, there 
exists only one $b > 0$ such that $\Phi_{k}^{>}(a,b)(T_{1})=T_{2}$, 
which we denote by $f_{k}^{>}(a)$. Applying the implicit function 
theorem together with lemma \ref{lem:2.3} to the relation 
$\Phi_{k}^{>}(a,b)(T_{1})=T_{2}$, one deduces that $I$ is open and 
that $f_{k}^{>}$ is of class $C^{1}$, with $(f_{k}^{>})' < 0$. 
Moreover, using the monotonicity and continuity properties of 
$\Phi_{k}^{>}$, one easily verifies that $I$ is an interval. It 
remains to prove that $I$ is unbounded from above., i. e.  of the 
form $]\nu_{k}^{>}, + \infty[$, and that $f_{k}^{>}(x) \to + \infty$ 
as $x \to \nu_{k}^{>}$. To verify the unboundedness one observes that 
if $a \in I$ and $\bar{a} > a$, then
$$
\Phi_{k}^{>}(\bar{a},f_{k}^{>}(a))(T_{1}) < 
\Phi_{k}^{>}(a,f_{k}^{>}(a))(T_{1})=T_{2};
$$
since, by lemma \ref{lem:2.5}, $\Phi_{k}^{>} (\bar{a},b)$ is $> 
T_{2}$ for some (small) $b > 0$, the conclusion that $\bar{a} \in I$ 
follows by continuity. Now if $f_{k}^{>}(x) \to L < +\infty$ as $x \to 
\nu_{k}^{>}$, then, by continuity, we get $\Phi_{k}^{>} 
(\nu_{k}^{>},L)(T_{1})=T_{2}$, so that $\nu_k	^>	 \in I$, which contradicts the fact that $I$ is open. Finally the
inequalities $\nu_{k}^{>}
\geq 
\lambda_{1}^{m}$ and $\mu_{k}^{>} \geq \lambda_{1}^{n}$ follow from 
the general properties of $\Sigma$ mentioned at the beginning of 
section 3. Q. E. D. 

\vspace{5mm}

The next proposition gives some information on how these hyperbolic 
like curves are situated one with respect to the other.

\begin{prop}\label{prop:3.5}
    Let $k \geq 2$. Assume that $C_{k+1}^{>++}$ (or $C_{k+1}^{<++}$)is 
    nonempty. Then $C_{k}^{>++}$ and $C_{k}^{<++}$ are both nonempty. 
    Moreover the curves $f_{k}^{>}$ and $f_{k}^{<}$ are both strictly 
    below the curve $f_{k+1}^{>}$ (or $f_{k+1}^{<}$).
\end{prop}

{\bf Proof.} Let us assume $C_{k+1}^{>++}$ nonempty (the argument is 
similar if it is $C_{k+1}^{<++}$ which is nonempty). By lemma 
\ref{lem:3.3}, this hypothesis means 
\begin{equation}
    \lim_{a,b \to + \infty} \Phi_{k+1}^{>}(a,b)(T_{1}) < T_{2}.
    \label{eq:3.3}
\end{equation}
Let us suppose for instance $k$ even. Then 
$$
\Phi_{k+1}^{>}(a,b)(T_{1})=\varphi_{a}^{m}[\Phi_{k}^{>}(a,b)(T_{1})] > 
\Phi_{k}^{>}(a,b)(T_{1})
$$
and consequently, by (\ref{eq:3.3}), we get
$$
\lim_{a,b \to + \infty} \Phi_{k}^{>}(a,b)(T_{1}) < T_{2},
$$
which implies, by lemma \ref{lem:3.3}, $C_{k}^{>++} \neq 0$.
Similar argument for $k$ odd. Let us now prove that $C_{k}^{<++}$ is 
also non empty. The argument in fact is similar to the preceding one 
and is now based on the relations
$$
\Phi_{k+1}^{>}(a,b)(T_{1})=\Phi_{k}^{<}(a,b)[\varphi_{a}^{m}(T_{1})] > 
\Phi_{k}^{<}(a,b)(T_{1}).
$$
It remains to see that $f_{k}^{>}$ and $f_{k}^{<}$ both lie strictly 
below $f_{k+1}^{>}$. Suppose by contradiction that $f_{k}^{>}(x) \geq 
f_{k+1}^{>}(x)$ for some $x > \max(\nu_{k}^{>}, \nu_{k+1}^{>})$. Then
\begin{eqnarray*}
T_{2} &=& \Phi_{k+1}^{>} (x,f_{k+1}^{>}(x))(T_{1})
\geq \Phi_{k+1}^{>}(x,f_{k}^{>}(x))(T_{1})\\
&>& \Phi_{k}^{>}(x,f_{k}^{>}(x))(T_{1})=T_{2},
\end{eqnarray*}
a contradiction. Similar argument for $f_{k}^{<}$. Q. E. D.

\vspace{5mm}

The rest of this section is mainly concerned with some results on the 
number of these hyperbolic like curves. Theorem 3.6 below is the 
converse of the first part of proposition 3.5.

\begin{theo}\label{theo:3.6}
Suppose $m^{+}(t) \not \equiv 0$ and $n^{+}(t) \not \equiv 0$ on 
$]T_{1},T_{2}[$. Then $C_{2}^{>++}$ or $C_{2}^{<++}$ is nonempty. 
Moreover if for some $k \geq 2$, $C_{k}^{>++}$ and $C_{k}^{<++}$ are 
both nonempty, then $C_{k+1}^{>++}$ or $C_{k+1}^{<++}$ is nonempty.
\end{theo}

Combining with remark 3.2, we get the following
\begin{coro}\label{coro:3.7}
If both $m(t)$ and $n(t)$ change sign, then $\Sigma^{*}$ contains at 
least one hyperbolic like curve in each quadrant.
\end{coro}

\begin{coro}\label{coro:3.8}
Suppose $m^{+}(t) \not \equiv 0$ and $n^{+}(t) \not \equiv 0$. Then 
either (i) $C_{k}^{>++}$ and $C_{k}^{<++}$ are nonempty for all $k$, 
or (ii) for some $k_{0} \geq 1$, $C_{k}^{>++}$ and $C_{k}^{<++}$ are 
nonempty for all $k \leq k_{0}$, $C_{k_0+1}^{>++}$ or $C_{k_0+1}^{<++}$ 
is nonempty, the other one being empty, and $C_{k}^{>++}$ and 
$C_{k}^{<++}$ are empty for all $k \geq k_{0}+2$. 
\end{coro}

So, if $m(t)$ and $n(t)$ both change sign, then each quadrant contains a 
(non zero) odd or infinite number of nonempty sets $C_{k}^{>}$, 
$C_{k}^{<}$ with $k \geq 2$. We will see at the end of this section 
that all cases can effectively happen (cf. proposition 3.12).

\begin{rem}\label{rem:3.9}
{\rm It may happen that $C_{k}^{>++}=C_{k}^{<++}$. For instance on 
$[-1,+1]$, if $p(-t)=p(t)$, $q(-t)=q(t)$, $m(-t)=m(t)$ and 
$n(-t)=n(t)$, then $C_{k}^{>++}=C_{k}^{<++}$ for all $k$ even. On the 
other hand, in the one weight problem on $[0,2\pi]$, if $m(t)=\sin t$ 
on $[0,\pi]$ and $m(t)=0$ on $[\pi,2\pi]$, then $C_{2}^{>++} \neq 
C_{2}^{<++}$. In fact theorem \ref{theo:5.1} implies that 
$C_{2}^{>++}$ and $C_{2}^{<++}$ have different asymptotic behaviours.}
\end{rem}

{\bf Proof of Theorem \ref{theo:3.6}.} Let us first consider 
$C_{2}^{>++}$ and $C_{2}^{<++}$. We will use the notation 
(\ref{eq:2.6}) as well as a similar one for the weight $n(t)$~: 
$\beta_{s}^{>} := \inf \{t > s : n(t) > 0\}$. Since $m^{+}(t) \not 
\equiv 0$ and $n^{+}(t) \not \equiv 0$, we have $\alpha_{T_{1}}^{>} < 
T_{2}$ and $\beta_{T_{1}}^{>} < T_{2}$. Two cases are distinguished~: 
(i) $\alpha_{T_{1}}^{>} \leq \beta_{T_{1}}^{>}$ or (ii) 
$\beta_{T_{1}}^{>} < \alpha_{T_{1}}^{>}$. Consider case(i). By the 
definition of $\beta_{T_{1}}^{>}$, 
\begin{equation}
    n^{+}(t) \not \equiv 0 \mbox{ on } ]\beta_{T_{1}}^{>}, 
    \beta_{T_{1}}^{>}+\epsilon[
    \label{eq:3.4}
\end{equation}
for any $\epsilon > 0$. Since, by lemma \ref{lem:2.4}, 
$\varphi_{a}^{m}(T_{1})$ converges to $\alpha_{T_{1}}^{>} \leq 
\beta_{T_{1}}^{>}$ as $a \to + \infty$, we deduce from (\ref{eq:3.4}) 
and lemma \ref{lem:2.4} that 
$$
\lim_{a,b \to + \infty}\varphi_{b}^{n}(\varphi_{a}^{m}(T_{1})) \leq 
\beta_{T_{1}}^{>} < T_{2}.
$$
Lemma \ref{lem:3.3} then implies $C_{2}^{>++} \neq \emptyset$. In 
case (ii), a similar argument yields $C_{2}^{<++} \neq \emptyset$.

Let us now turn to the study of $C_{k+1}^{>++}$ and $C_{k+1}^{<++}$, 
assuming $k \geq 2$, $C_{k}^{>++} \neq \emptyset$ and $C_{k}^{<++} 
\neq \emptyset$. This hypothesis means, by lemma \ref{lem:3.3}, that 
$$
\alpha : = \lim_{a,b \to + \infty} \Phi_{k}^{>}(a,b)(T_{1}) < T_{2} 
\mbox{ and } \beta := \lim_{a,b \to + \infty} \Phi_{k}^{<} (a,b) 
(T_{1}) < T_{2}.
$$
We again distinguish two cases : (i) $\alpha \leq \beta$ or (ii) 
$\beta < \alpha$. Consider case (i) and let us first assume $k$ odd. 
So 
$$
\Phi_{k}^{<}(a,b)(T_{1})=\varphi_{b}^{n}\left[ 
\Phi_{k-1}^{<}(a,b)(T_{1})\right] \to \beta \mbox{ as } a,b \to + 
\infty.
$$
This implies that
\begin{equation}
    n^{+}(t) \not \equiv 0 \mbox{ on } ]\beta, \beta+\epsilon[
    \label{eq:3.5}
\end{equation}
for any $\epsilon > 0$. Indeed, for $a,b$ sufficiently large, one 
has $\beta < \varphi_{b}^{n}[\Phi_{k-1}^{<} (a,b)(T_{1})] < \beta + 
\epsilon$, which implies $n^{+}(t) \not \equiv 0$ on 
$[\Phi_{k-1}^{<}(a,b)(T_{1}),\beta+\epsilon]$, and (\ref{eq:3.5}) 
follows. Now, $\Phi_{k}^{>}(a,b)(T_{1})$ converges to $\alpha \leq 
\beta$ as $a,b \to +\infty$, and so we deduce from (\ref{eq:3.5}) and 
lemma \ref{lem:2.4} that 
$$
\lim_{a,b \to +\infty}\varphi_{b}^{n}[\Phi_{k}^{>}(a,b)(T_{1})]\leq 
\beta < T_{2}.
$$
Consequently, by lemma \ref{lem:3.3}, $C_{k+1}^{>++} \neq 0$. A 
similar argument for $k$ even would also lead to $C_{k+1}^{>++} \neq 
\emptyset$. And in case (ii) one would get $C_{k+1}^{<++} \neq 
\emptyset$. Q. E. D.

\vspace{3mm}

We will now give a sufficient (and almost necessary) condition on the 
weights $m(t)$ and $n(t)$ in order that $\Sigma^{*} \cap (\dR^{+} \times \dR^{+})$ contains an 
infinite number of hyperbolic like curves.

\begin{theo}
\label{theo:3.10}
If $m^{+}(t) \cdot n^{+}(t) \not \equiv 0$ on $]T_{1},T_{2}[$, then 
$C_{k}^{>++}$ and $C^{<++}_{k}$ are nonempty for all $k$. Conversely 
suppose that the set where $m(t)$ (or $n(t))$ is $>0$ is made of a 
finite union of intervals. Under this hypothesis, if $C_{k}^{>++}$ 
and $C_{k}^{<++}$ are nonempty for all $k$, then $m^{+}(t) \cdot 
n^+(t) \not \equiv 0$ on $]T_{1},T_{2}[$.
\end{theo}

{\bf Proof.} If $m^{+}(t) \cdot n^{+}(t) \not \equiv 0$, then one can 
find $T_{1} \leq t_{1} < t_{2} \leq T_{2}$ such that $m(t)$ and $n(t)$ 
are both $>0$ on $[t_{1}, t_{2}]$. Lemma \ref{lem:2.4} then implies 
that for each $k$, $\Phi_{k}^{>}(a,b)(t_{1})$ and 
$\Phi_{k}^{<}(a,b)(t_{1}) \to t_{1}$ as $a,b \to + \infty$. By 
monotonicity we deduce that for each $k$, the limits of 
$\Phi_{k}^{>}(a,b)(T_{1})$ and $\Phi_{k}^{<}(a,b)(T_{1})$ are $\leq 
t_{1} < T_{2}$ as $a,b \to + \infty$. Lemma \ref{lem:3.3} then yields 
the conclusion.

Conversely let us assume that $C_{k}^{>++}$ and $C_{k}^{<++}$ are 
nonempty for all $k$ and that, say, $m(t)$ is $>0$ on $I_{1} \cup 
I_{2} \cup \ldots \cup I_{r}$ and $\leq 0$ outside this set, where 
$I_{i}$ is an interval of extremities $s_{i}, t_{i}$, with 
$$
T_{1} \leq s_{1} < t_{1} \leq s_{2} < t_{2} \leq \ldots \leq s_{r} 
<t_{r} \leq T_{2}.
$$
Suppose by contradiction that $m^{+}(t) \cdot n^{+}(t) \equiv 0$ on 
$]T_{1}, T_{2}[$. So $n(t) \leq 0$ on each $[s_{i},t_{i}]$. This 
implies that 
$$
\varphi_{a}^{m}(T_{1}) \geq s_{1}, \; \varphi_{b}^{n}(s_{i}) \geq 
t_{i}, \, \varphi_{a}^{m}(t_{i}) \geq s_{i+1}, \, 
\varphi_{a}^{m}(t_{r})\geq T_{2}
$$
for all $a,b \geq 0$ and each $i=1,\ldots,r$. Consequently 
$\Phi_{k}^{>}(a,b)(T_{1}) \geq T_{2}$ for $k \geq 2r+2$ and for all 
$a,b \geq 0$. Lemma \ref{lem:3.3} then implies $C_{k}^{>++} = 
\emptyset$ for $k \geq 2r+2$, a contradiction. Similar argument if it 
is $n(t)$ which is $>0$ on a finite union of intervals. Q. E. D.

\begin{rem}
\label{rem:3.11}
{\rm Some hypothesis on the bumps of the weights is needed in the 
second part of theorem \ref{theo:3.10}. For example if $[T_{1}, 
T_{2}]=[0,\pi]$ and $m(t)=-n(t)=t \sin (1/t)$, then $m^{+}(t)\cdot 
n^{+}(t) \equiv 0$ although $C_{k}^{>++}$ and $C_{k}^{<++}$ are 
nonempty for all $k$. This latter fact follows from lemma \ref{lem:3.3} 
since here, $\Phi_{k}^{>}(a,b)(0)$ and $\Phi_{k}^{<}(a,b)(0) \to 0$ 
as $a,b \to +\infty$. More general hypothesis on the weights are 
considered in \cite{Al1} which guarantee the validity of the second 
part of theorem \ref{theo:3.10}.

\vspace{3mm}

To conclude this section, we show that all cases can effectively 
happen with respect to the numbers of nonempty sets $C_{k}^{>}, 
C_{k}^{<}$ with $k \geq 2$ in the various quadrants.
}
\end{rem}

\begin{prop}
    \label{prop:3.12}
Given $K,L,M,N$ in $\{0,1,2,\ldots,+\infty\}$, there exist weights 
$m(t)$, $n(t)$, which both change sign, such that the total number of 
nonempty sets $C_{k}^{>++}$ and $C_{k}^{<++}$ (resp. $C_{k}^{>--}$ 
and $C_{k}^{<--}$, $C_{k}^{>+-}$ and $C_{k}^{<+-}$, $C_{k}^{>-+}$ and 
$C_{k}^{<-+}$) with $k \geq 2$ is equal to $2K+1$ (resp. $2L+1, 2M+1, 
2N+1)$.
\end{prop}

{\bf Proof.} The proof is based on the observation that the number
of nonempty sets $C_{k}^{>}$ and $C_{k}^{<}$ with $k \geq 2$ in 
$\dR^{+} \times \dR^{+}$ (resp. $\dR^{-} \times \dR^{-}$, $\dR^{+} 
\times \dR^{-}$, $\dR^{-} \times \dR^{+}$) only depends on the 
relative position of the positive (resp. negative, positive, 
negative) bumps of $m$ and of the positive (resp. negative, negative, 
positive) bumps of $n$. A positive bump of $m$ which intersects a 
positive bump of $n$ leads, by theorem \ref{theo:3.10}, to an 
infinity of nonempty sets $C_{k}^{>}$ and $C_{k}^{<}$ with $k \geq 2$ 
in $\dR^{+} \times \dR^{+}$. A succession of $K+2$ alternating 
positive bumps of $m$ and positive bumps of $n$ (i.e., for instance, 
$m(t)=(\sin t)^{+}$ and $n(t)=(\sin(t+\pi))^{+}$ on 
$[0,(K+2)\pi]$) leads to exactly $2K+1$ nonempty sets $C_{k}^{>}$ 
and $C_{k}^{<}$ with $k \geq 2$ in $\dR^{+} \times \dR^{+}$ (a simple 
consequence of lemmas \ref{lem:2.4} and \ref{lem:3.3}). This latter 
construction of alternating bumps will be used repeatedly below.

To carry through the details of the proof, we will distinguish 
several cases according to the number of quadrants having to contain 
an infinite number of sets $C_{k}^{>}$ and $C_{k}^{<}$ with $k \geq 2$.

\vspace{3mm}

{\sl Case 1 :} four ``infinite'' quadrants, i.e. $K=L=M=N=+\infty$. 
Applying theorem \ref{theo:3.10} and remark \ref{rem:3.2}, one easily 
verifies that this situation occurs for instance on $[0,2\pi]$ if we 
take $m(t)=\sin t$ and $n(t)=\sin (t+\pi/2)$.

\vspace{3mm}

{\sl Case 2 :} one ``finite'' quadrant, three ``infinite'' quadrants. 
Applying remark \ref{rem:3.2}, one sees that it suffices to construct 
an example with $K$ given in $\{0,1,2,\ldots\}$ and $L=M=N=+\infty$. 
Using theorem \ref{theo:3.10} and lemmas \ref{lem:2.4} and 
\ref{lem:3.3}, one easily verifies that this situation occurs for 
instance on $[0,(K+2)\pi]$ if we take $m(t)=\sin t$ on $[0,2\pi]$, 
$n(t)=-\sin 2t/3$ on $[0,3\pi/2]$, $n(t)=-\sin 2(t-3\pi/2)$ on 
$[3\pi/2,2\pi]$, and then, on $[2\pi,(K+2)\pi],K$ alternating 
positive bumps for $m$ and positive bumps for $n$ (starting with $m$).

\vspace{3mm}

{\sl Case 3 :}
two ``finite'' quadrants, two ``infinite'' quadrants. 
Applying remark \ref{rem:3.2}, one sees that it suffices to construct 
an example with $K,L$ given in $\{0,1,2,\ldots\}$ and $M=N=+\infty$, 
an example with $K,M$ given in $\{0,1,2,\ldots\}$ and $L=N=+\infty$, 
and an example with $K,N$ given in $\{0,1,2,\ldots\}$ and 
$L=M=+\infty$. We will describe below the construction relative to 
$K,L$ finite (and $K \geq L$, $L$ even). The other cases can be 
treated along similar lines.  The desired situation occurs for 
instance on $[0,(K+2)\pi)$ if we take $m(t)=\sin t$ on $[0,(L+2)\pi]$, 
$n(t)=-\sin t$ on $[0,(L+2)\pi]$, and then, on 
$[(L+2)\pi,(K+2)\pi]$, $K-L$ alternating positive bumps for $m$ and 
positive bumps for $n$ (starting with $m$).

\vspace{3mm}

{\sl Case 4 :} three ``finite'' quadrants, one ``infinite'' quadrant. 
Applying remark \ref{rem:3.2}, one sees that it suffices to construct 
an example with $K=+\infty$ and $L,M,N$ given in $\{0,1,2,\ldots\}$. 
We will describe below the contsruction for $L,M,N$ even $\neq 0$. The 
other cases can be treated along similar lines. The desired situation 
occurs for instance if we take $m(t)=n(t)=\sin t$ on $[0,\pi]$, and 
then $L$ alternating negative bumps for $m$ and $n$ (starting with 
$m$), followed by $N$ alternating positive bumps for $n$ and negative 
bumps for $m$ (starting with $n$), followed by $M$ alternating 
positive bumps for $m$ and negative bumps for $n$ (starting with 
$m$).

\vspace{3mm}

{\sl Case 5 :} four ``finite'' quadrants, i.e. $K, L, M, N$ given 
in $\{0,1,2,\ldots\}$. We will describe below the construction for 
$K,L,M,N$ even $\neq 0$. The other cases can be treated along similar 
lines. The desired situation occurs for instance if we take $K$ 
alternating positive bumps for $m$ and $n$ (starting with $m$), 
followed by $M$ alternating positive bumps for $m$ and negative bumps 
for $n$ (starting with $m$), followed by one positive bump for $m$, 
followed by $L$ alternating negative bumps for $m$ and $n$ (starting 
with $m$), followed finally by $N$ alternating positive bumps for $n$ 
and negative bumps for $m$ (starting with $n$). Q. E. D.

\begin{exam}\label{exam:3.13}
{\rm If we take a positive bump for $m$, followed by a positive bump 
for $n$, followed by a negative bump for $m$, followed by a negative 
bump for $n$, (i.e., for instance, $m(t)=(\sin t)^+$, $n(t)=(\sin t)^-$ on$ [0,2\pi]$, and 
$m(t)=-(\sin t)^+$,
$n(t)=-(\sin t)^-$ on
$[2\pi,4\pi]$),  then
$\Sigma^{*}$ contains exactly one hyperbolic like  curve in each quadrant.
}
\end{exam}

Proposition \ref{prop:3.12} concerns weights which both change sign. 
Similar arguments lead to the following two propositions.

\begin{prop}\label{prop:3.14}
Given a pair $(A,B)$ of quadrants different from the diagonal 
pairs $(\dR^{+} \times \dR^{+}, \dR^{-} \times \dR^{-})$ or $(\dR^{+} 
\times \dR^{-}, \dR^{-} \times \dR^{+})$, and given $R,S \in 
\{0,1,2,\ldots,+\infty\}$, there exist weights $m(t)$ and $n(t)$, one 
of them which changes sign and the other one which does not, such that the 
total number of nonempty sets $C_{k}^{>}$ and $C_{k}^{<}$ with $k \geq 
2$ in the quadrant $A$ (resp. $B$) is equal to $2R+1$ (resp. $2S+1)$, 
while no such sets appear in the other two quadrants.
\end{prop}

\begin{prop}\label{prop:3.15}
Given a quadrant and $R \in \{0,1,2,\ldots,+\infty\}$, there exist 
weights $m(t)$ and $n(t)$, which do no change sign, such that the 
total number of nonempty sets $C_{k}^{>}$ and $C_{k}^{<}$ with $k 
\geq 2$ in the given quadrant is $2R+1$, while no such sets appear in 
the other three quadrants.
\end{prop}

\section{The one weight problem}
\setcounter{equation}{0}

We now consider problem (\ref{eq:1.1}) with $L$ as in the 
introduction and $m \in C[T_{1},T_{2}]$, $m(t) \not \equiv 0$.

In this case the Fu\v cik spectrum $\Sigma$ is symmetric with respect 
to the diagonal $a=b$ in the $(a,b)$ plane. In fact replacing $u$ by 
$-u$ in (\ref{eq:1.1}) shows that $C_{k}^{>++}$ (resp. $C_{k}^{>--}, 
C_{k}^{>+-}, C_{k}^{>-+}$) is the symmetric of $C_{k}^{<++}$ (resp. 
$C_{k}^{<--}, C_{k}^{<-+}, C_{k}^{<+-}$). Moreover 
$(\lambda_{k}^{m}, \lambda_{k}^{m})\in C_{k}^{>++} \cap C_{k}^{<++}$ 
and $(\lambda_{-k}^{m}, \lambda_{-k}^{m}) \in C_{k}^{>--} \cap 
C_{k}^{<--}$ for all $k$ (depending of course whether $m$ is $\geq 
0$, $\leq 0$ or changes sign); in particular the corresponding sets 
$C_{k}^{>}$ and $C_{k}^{<}$ in $\dR^{+} \times \dR^{+}$ and - or 
$\dR^{-} \times \dR^{-}$ are nonempty.
As in section 3, all the sets $C_{k}^{>++}, C_{k}^{<++}, C_{k}^{>--}, 
C_{k}^{<--}, C_{k}^{>+-},C_{k}^{<+-}$,\\$C_{k}^{>-+}, C_{k}^{<-+}$ with 
$k \geq 2$, when nonempty, are hyperbolic like curves (cf. theorem 
3.4).

The case where the weight $m(t)$ does not change sign is simpler~: if 
$m \geqs{\not \equiv} 0$ (resp. $m \leqs{\not \equiv} 0$), then 
$\Sigma^{*}$ is made of the infinite sequence of hyperbolic like 
curves $C_{k}^{>++}$ and $C_{k}^{<++}$ (resp. $C_{k}^{>--}$ and 
$C_{k}^{<--}$), $k=2,3,\ldots$

From now on in this section we will assume that 
\begin{equation}
    m(t) \mbox{ changes sign in } ]T_{1},T_{2}[.
    \label{eq:4.1}
\end{equation}
Under this assumption, $\Sigma^{*}$ contains the infinite sequence of 
hyperbolic like curves \\$C_{k}^{>++}, C_{k}^{<++}, C_{k}^{>--}, 
C_{k}^{<--}$, $k=2,3,\ldots$. Moreover, by remark \ref{rem:3.2} and 
theorem \ref{theo:3.6}, hyperbolic like curves also appear in 
$\dR^{+} \times \dR^{-}$ and $\dR^{-} \times \dR^{+}$. These 
additional curves in $\dR^{+} \times \dR^{-}$ are symmetric of those in $\dR^-	 \times \dR^+	$, 
and their distribution 
is given by corollary \ref{coro:3.8}~: either (i) $C_{k}^{>+-}$ and 
$C_{k}^{<+-}$ are nonempty for all $k$, or (ii) for some $k_{0} \geq 
1$, $C_{k}^{>+-}$ and $C_{k}^{<+-}$ are nonempty for all $k \leq 
k_{0}$, $C_{k_{0}+1}^{>+-}$ or $C_{k_{0}+1}^{<+-}$ is nonempty, the 
other one being empty, and $C_{k}^{>+-}$ and $C_{k}^{<+-}$ are empty 
for all $k \geq k_{0} +2$.

The rest of this section is concerned with showing that the number of 
these additional curves in $\dR^{+}\times \dR^{-}$ and $\dR^-	 \times \dR^+	$ 
is directly related 
to the number of changes of sign of the weight $m(t)$. First we have 
to make precise this notion of {\sl number of changes of sign.}

\begin{defi}\label{defi:4.1}
   {\rm  Let $s \in ]T_{1},T_{2}[$. We say that $s$ is a {\sl simple} point 
    of change of sign of $m$ if there exists $T_{1} < s' \leq s$ and 
    $\epsilon_{0} > 0$ such that either (i) $m \leqs{\not \equiv} 0$ 
    on $]s'-\epsilon,s'[$ for all $0 < \epsilon < \epsilon_{0}$, $m 
    \equiv 0$ on $[s',s]$ and $m \geqs{\not \equiv} 0$ on $]s,s+\epsilon[$ for all 
$0 < \epsilon < \epsilon_0$, or (ii) $m \geqs{\not \equiv} 0$ on
    $]s'-\epsilon,s'[$ for all $0 <  \epsilon < \epsilon_{0}$, 
    $m\equiv 0$ on $[s',s]$ and $m \leqs{\not \equiv} 0$ on 
    $]s,s+\epsilon[$ for all $0 < \epsilon < \epsilon_{0}$.}
\end{defi}

\begin{defi}\label{defi:4.2}
    {\rm Let $s \in [T_{1}, T_{2}]$. We say that $s$ is a {\sl multiple} 
    point of change of sign of $m$ if either (i) $s > T_1$ and $m^{+}$ and $m^{-}$ 
    are $\not \equiv 0$ on $]s-\epsilon, s[ \, \cap \, [T_{1}, T_{2}]$ for 
    any $\epsilon > 0$, or (ii) $s < T_2$ and $m^{+}$ and $m^{-}$ are $\not \equiv 
    0$ on $]s,s+\epsilon[ \, \cap \, [T_{1}, T_{2}]$ for any $\epsilon > 0$.}
\end{defi}

\begin{defi}\label{defi:4.3}
    {\rm If $m$ admits at least one multiple point of change of sign, then 
    we say that {\sl the number of changes of sign of $m$ is} 
    $+\infty$. If $m$ does not admit any multiple point of change of 
    sign, then it is easily verified that, under (\ref{eq:4.1}), $m$ 
    admits a nonzero finite number $N$ of simple points of change of 
    sign. In this case we say that {\sl the number of changes of sign 
    of $m$ is} $N$.}
\end{defi}

\begin{theo}\label{theo:4.4}
    Assume (\ref{eq:4.1}). Let $N \in \{1,2,\ldots,+\infty\}$ be the 
    number of changes of sign of $m$. Then the total number of 
    nonempty sets $C_{k}^{>+-}$ and $C_{k}^{<+-}$ with $k \geq 2$ is 
    equal to $2N-1$.
\end{theo}

{\bf Proof.} Consider first the situation where $N=+\infty$ and let $s 
\in [T_{1}, T_{2}]$ be a multiple point of change of sign of $m$. If 
$s < T_{2}$, then one easily verifies, using lemma \ref{lem:2.4}, that 
for any $k=2,3,\ldots,$
$$
\biindice{\lim}{a \to + \infty}{b \to - \infty} 
\Phi_{k}^{>}(a,b)(T_{1}) \leq s \mbox{ and } \biindice{\lim}{a \to + 
\infty}{b \to -\infty} \Phi_{k}^{<} (a,b) (T_{1}) \leq s.
$$
If $s=T_{2}$, then one gets that the above limits are $< s$. So, in 
any case, lemma \ref{lem:3.3} implies that $C_{k}^{>+-}$ and 
$C_{k}^{<+-}$ are nonempty for all $k$. Consider now the situation 
where $N$ is finite. Then there exist 
$$
s_{1} = T_{1} < s_{2} < \ldots < s_{N+1} < s_{N+2} = T_{2}
$$
such that either $m$ is $\geqs{\not \equiv} 0$  on $[s_{1}, s_{2}]$, $\leqs{\not \equiv} 0$ 
on $[s_{2},s_{3}]$, $\geqs{\not \equiv}0$ on $[s_{3},s_{4}]$,
$\ldots$, or $m$ is $\leqs{\not \equiv} 0$ $[s_{1},s_{2}]$, 
$\geqs{\not \equiv} 0$ on $[s_{2},s_{3}]$, $\leqs{\not \equiv} 0$ on 
$[s_{3}, s_{4}]$, $\ldots$ Let us deal with the first case with, say, 
$N$ even, so that $m$ is $\geqs{\not \equiv} 0$ on $[s_{N+1}, 
s_{N+2}]$ (the other cases  are treated similarly). Using lemma 
\ref{lem:2.4}, one gets that 
$$
\biindice{\lim}{a \to +\infty}{b \to -\infty} 
\Phi_{N+1}^{>}(a,b)(T_{1}) < T_{2}, \biindice{\lim}{a \to + 
\infty}{b \to - \infty} \Phi_{N+2}^{>}(a,b)(T_{1}) = T_{2},
$$
$$
\biindice{\lim}{a \to +\infty}{b \to -\infty} 
\Phi_{N}^{<}(a,b)(T_{1}) < T_{2}, \biindice{\lim}{a \to + 
\infty}{b \to - \infty} \Phi_{N+1}^{<}(a,b)(T_{1}) = T_{2}.
$$
Consequently, lemma \ref{lem:3.3} implies that $C_{k}^{>+-}$ and 
$C_{k}^{<+-}$ are nonempty for $k=2,\ldots,N$, that $C_{N+1}^{>+-}$ 
is nonempty, that $C_{N+1}^{<+-}$ is empty, and that all $C_{k}^{>+-}$ and $C_k^{<+-}$
are empty for $k \geq N+2$. The conclusion follows. Q. E. D.

\vspace{3mm}

Theorem \ref{theo:4.4} implies a result analogous to proposition 
\ref{prop:3.12}~: given $N \in \{1,2,\ldots,+\infty\}$, there exists a 
weight $m(t)$ which changes sign such that $\dR^{+} \times \dR^{-}$ 
exactly contains $2N-1$ nonempty sets $C_{k}^{>+-}$ and $C_{k}^{<+-}$ 
with $k \geq 2$.

\begin{rem}\label{rem:4.5}
    {\rm It may happen that $C_{k}^{>+-}=C_{k}^{<+-}$ (and then of course, 
    by symmetry, $C_{k}^{>-+}=C_{k}^{<-+}$). This is the case for 
    instance for all $k$ even if $[T_{1},T_{2}]=[-1,+1]$, 
    $p(-t)=p(t)$, $q(-t)=q(t)$ and $m(-t)=m(t)$. In fact, in this 
    example, one also has $C_{k}^{>++}=C_{k}^{<++}$ and 
    $C_{k}^{>--}=C_{k}^{<--}$ for all $k$ even. On the other hand, 
    if $[T_{1}, T_{2}]=[0,2\pi]$ and $m(t)=\sin t$, then 
    $C_{2}^{>+-} \neq C_{2}^{<+-}$. In fact, theorem \ref{theo:5.1} 
    implies that these two curves have different asymptotic behaviours.}
\end{rem}

\section{Asymptotic behaviour of the first curves}
\setcounter{equation}{0}

In this section we return to the two weights problem (\ref{eq:1.3}) with $L$ as in the 
introduction and $m,n \in C[T_{1},T_{2}]$, $m(t)$ and $n(t) \not 
\equiv 0$. Our purpose is to investigate the asymptotic behaviour of 
the first hyperbolic like curves $C_{2}^{>++}, 
C_{2}^{<++},C_{2}^{>--}, C_{2}^{<--}, C_{2}^{>+-}, C_{2}^{<+-}, 
C_{2}^{>-+}, C_{2}^{<-+}$. We recall that at least one such curve 
appears in $\dR^{+} \times \dR^{+}$ (resp. $\dR^{-} \times \dR^{-}, 
\dR^{+} \times \dR^{-}, \dR^{-} \times \dR^{+})$ if (and only if) 
$m^{+} \not \equiv 0$ and $n^{+} \not \equiv 0$ (resp. $m^{-} \not 
\equiv 0$ and $n^{-} \not \equiv 0$, $m^{+} \not \equiv 0$ and $n^{-} 
\not \equiv 0$, $m^{-} \not \equiv 0$ and $n^{+} \not \equiv 0$). In 
particular, in the one weight problem (\ref{eq:1.1}), at least one 
such curve appears in each quadrant if $m$ changes sign.

Some notations are needed to state our result. Let us define 
\begin{eqnarray*}
    T_{1}^{m>} & := & \inf \{ t \in ]T_{1}, T_{2}[ : m(t) > 0\},  \\
    T_{2}^{m>} & := & \sup \{ t \in ]T_{1}, T_{2}[ : m(t) > 0\},  \\
    T_{1}^{m<} & := & \inf \{ t \in ]T_{1}, T_{2}[ : m(t) < 0\},  \\
    T_{2}^{m<} & := & \sup \{ t \in ]T_{1}, T_{2}[ : m(t) < 0\},
\end{eqnarray*}
and similarly for $n : T_{1}^{n>}, T_{2}^{n>}, T_{1}^{n<}, T_{2}^{n<}$. 
Note that some of these quantities may not have sense in $[T_1,T_2]$~: for instance 
$T_{1}^{m>}$ and $T_{2}^{m>}$ are only defined if $m^{+} \not \equiv 
0$, and similarly for the others. We will also denote by 
$\lambda_{1}(m,]t_{1},t_{2}[$) (resp. $\lambda_{-1}(m,]t_{1},t_{2}[$) 
the first positive (resp. negative) eigenvalue of $L$ with Dirichlet 
boundary condition on the interval $]t_{1}, t_{2}[$ for the weight $m$. 
Again here this quantity makes sense only if $m$ has a nontrivial 
positive (resp. negative) part on $]t_{1},t_{2}[$.

\begin{theo}\label{theo:5.1}
    \begin{enumerate}
        \item[(i)] If $C_{2}^{>++} \neq \emptyset$, then 
        $C_{2}^{>++}$ is asymptotic to the lines $\dR \times 
        \lambda_{1} (n, ]T_{1}^{m>}, T_{2}[$) and 
        $\lambda_{1}(m,]T_{1},T_{2}^{n>}[) \times \dR$.
    
        \item[(ii)] If $C_{2}^{<++} \neq 0$, then $C_{2}^{<++}$ is 
        asymptotic to the lines $\dR \times 
        \lambda_{1}(n,]T_{1},T_{2}^{m>}[)$ and 
        $\lambda_{1}(m,]T_{1}^{n>},T_{2}[)\times \dR$.
    
        \item[(iii)] If $C_{2}^{>--} \neq \emptyset$, then 
        $C_{2}^{>--}$ is asymptotic to the lines $\dR \times 
        \lambda_{-1}(n, ]T_{1}^{m<}, T_{2}[$) and 
        $\lambda_{-1}(m,]T_{1},T_{2}^{n<}[)\times \dR$.
    
        \item[(iv)] If $C_{2}^{<--} \neq \emptyset$, then 
        $C_{2}^{<--}$ is asymptotic to the lines $\dR \times 
        \lambda_{-1}(n, ]T_{1},T_{2}^{m<}[$) and $\lambda_{-1}(m, 
        ]T_{1}^{n<},T_{2}[) \times \dR$.
    
        \item[(v)] If $C_{2}^{>+-} \neq \emptyset$, then 
        $C_{2}^{>+-}$ is asymptotic to the lines $\dR \times 
        \lambda_{-1}(n,]T_{1}^{m>},T_{2}[)$ and $\lambda_{1}(m, 
        ]T_{1},T_{2}^{n<}[) \times \dR$.
    
        \item[(vi)] If $C_{2}^{<+-} \neq \emptyset$, then 
        $C_{2}^{<+-}$ is asymptotic to the lines $\dR \times 
        \lambda_{-1}(n,]T_{1},T_{2}^{m>}[$) and 
        $\lambda_{1}(m,]T_{1}^{n<},T_{2}[) \times \dR$.
    
        \item[(vii)] If $C_{2}^{>-+} \neq \emptyset$, then 
        $C_{2}^{>-+}$ is asymptotic to the lines $\dR \times 
        \lambda_{1} (n,]T_{1}^{m<},T_{2}[)$ and 
        $\lambda_{-1}(m,]T_{1},T_{2}^{n>}[) \times \dR$.
    
        \item[(viii)] If $C_{2}^{<-+} \neq \emptyset$, then 
        $C_{2}^{<-+}$ is asymptotic to the lines $\dR \times 
        \lambda_{1}(n,]T_{1},T_{2}^{m<}[$) and 
        $\lambda_{-1}(m,]T_{1}^{n>},T_{2}[) \times \dR$.
    \end{enumerate}
\end{theo}

\begin{rem}\label{rem:5.2}
    {\rm In (i) above, the hypothesis $C_{2}^{>++} \neq \emptyset$ 
    means that 
    \begin{equation}
        \varphi_{b}^{n}(\varphi_{a}^{m}(T_{1}))=T_{2} \mbox{ for some } 
        a,b > 0.
        \label{eq:5.1}
    \end{equation}
    This implies that $n^{+} \not \equiv 0$ on $]T_{1}^{m>}, T_{2}[$, 
    and consequently the eigenvalue 
    $\lambda_{1}(n,]T_{1}^{m>},T_{2}[$) which appears in (i) is 
    well-defined. (\ref{eq:5.1}) also implies that $m^{+} \not \equiv 
    0$ on $]T_{1},T_{2}^{n>}[$, and consequently the eigenvalue 
    $\lambda_{1}(m, ]T_{1},T_{2}^{n>}[$) which appears in (i) is 
    well-defined. Similar observation for the other curves.
    }
\end{rem}

{\bf Proof of Theorem \ref{theo:5.1}.} By using remark \ref{rem:3.2}, 
one sees that the statements relative to the curves in $\dR^{-} \times 
\dR^{-}$, $\dR^{+} \times \dR^{-}$ and $\dR^{-} \times \dR^{+}$ can 
be deduced from those relative to the curves in $\dR^{+} \times \dR^{+}$. 
Dealing with $\dR^{+} \times \dR^{+}$, we will only consider 
$C_{2}^{>++}$ ($C_{2}^{<++}$ can be treated similarly).

Let $(a,b) \in \dR^{+} \times \dR^{+}$. Then 
\begin{equation}
    (a,b) \in C_{2}^{>++} \Leftrightarrow 
    \varphi_{b}^{n}(\varphi_{a}^{m}(T_{1}))=T_{2} \Leftrightarrow 
    \varphi_{a}^{m}(T_{1})=(\varphi_{b}^{n})^{-1}(T_{2})
    \label{eq:5.2}
\end{equation}
If $a \to + \infty$ in (\ref{eq:5.2}), then, by lemma 
\ref{lem:2.4}, $\varphi_{a}^{m}(T_{1})\to T_{1}^{m>}$, and we deduce 
from part (ii) of lemma \ref{lem:5.3} below that $b \to \bar{b}$ 
with $\varphi_{\bar{b}}^{n} (T_{1}^{m>})=T_{2}$, i.e. 
$\bar{b}=\lambda_{1}(n,]T_{1}^{m>},T_{2}[)$. If $b \to + \infty$ in 
(\ref{eq:5.2}), then, by part (i) of lemma \ref{lem:5.3} below, 
$(\varphi_{b}^{n})^{-1}(T_{2}) \to T_{2}^{n>}$, and we deduce from 
part (ii) of lemma \ref{lem:5.3} below that $a \to \bar{a}$ 
with $\varphi_{\bar{a}}^{m}(T_{1})=T_{2}^{n>}$, i.e. $\bar{a} = 
\lambda_{1}(m,]T_{1},T_{2}^{n>}[$). Q. E. D.

\begin{lem}\label{lem:5.3}
    (i) $(\varphi_{a}^{m})^{-1}(T_{2})\to T_{2}^{m>}$ as $a \to 
    +\infty$. (ii) Let $\varphi_{a}^{m}(x)=y$ with $a > 0$. If $x \to 
    \bar{x}$ and $y \to \bar{y}$ with $\bar{x} \neq \bar{y}$, then 
    $a \to $ some $\bar{a} > 0$ which verifies 
    $\varphi_{\bar{a}}^{m}(\bar{x})=\bar{y}$.
\end{lem}

{\bf Proof.} Part (i) follows from lemma \ref{lem:2.4} by 
reversing the time. Part (ii) concerns the continuity of the 
function $(x,y) \to a > 0$ defined by the relation 
$\varphi_{a}^{m}(x)=y$. The implicit function theorem combined with 
(\ref{eq:2.5}) and an argument based on Rolle's theorem as in the 
proof of lemma \ref{lem:2.3} imply that this function is of class $C^{1}$. 
Q. E. D.

\vspace{3mm}

Part (ii) of lemma \ref{lem:5.3} can also be seen as a result on the 
continuous dependence of $\lambda_{1}(m,]x,y[)$ with respect to the 
domain $]x,y[$.

\begin{coro}\label{coro:5.4}
    If $m$ and $n$ both have compact support in $]T_{1}, T_{2}[$, then 
    none of the first hyperbolic like curves is asymptotic on any side 
    to the trivial horizontal and vertical  lines. The converse 
    implication also holds.
\end{coro}

{\bf Proof.} Assume that $m$ and $n$ have compact support and, to fix 
the ideas that both $m$ and $n$ change sign (similar  argument in 
the other cases). Then the quantities $T_{1}^{m>}, T_{1}^{m<}, 
T_{1}^{n>}, T_{1}^{n<}$ (resp. $T_{2}^{m>}, T_{2}^{m<}, T_{2}^{n>}, 
T_{2}^{n<}$) are all well-defined and $> T_{1}$ (resp. $< T_{2}$). It 
follows that the intervals like $]T_{1}^{m>}, T_{2}[, \ldots$ which 
appear in the statement of theorem \ref{theo:5.1} are all strictly 
smaller than $]T_{1}, T_{2}[$. By the strict monotonicity dependence 
of the first eigenvalues with respect to the domain, we deduce that 
the positive (resp. negative) first eigenvalue which appear in the 
statement of theorem \ref{theo:5.1} are all $>$ (resp. $<$) than the 
corresponding positive (resp. negative) first eigenvalues on 
$]T_{1},T_{2}[$. The conclusion of Corollary \ref{coro:5.4} follows. 
Consider now the converse implication and suppose, by contradiction, 
that, for instance, the support of $m$ hits $T_{1}$. Then the support 
of $m^{+}$ or the support of $m^{-}$ hits $T_{1}$. Suppose it is that 
of $m^{+}$. This implies that $T_{1}^{m>}=T_{1}$ and also that at 
least $C_{2}^{>++}$ or $C_{2}^{>+-}$ is $\neq \emptyset$. Suppose it 
is $C_{2}^{>++}$. 
Part (i) of theorem \ref{theo:5.1} then implies that $C_{2}^{>++}$ is 
asymptotic to $\dR \times \lambda_{1}(n, ]T_{1}, T_{2}[)=\dR \times 
\lambda_{1}^{n}$, a contradiction. Similar argument in the other 
cases. Q. E. D. 

\begin{coro}\label{coro:5.5}
    Consider the one weight problem (\ref{eq:1.1}) with $m \geqs{\not 
    \equiv} 0$. Then there exists $\epsilon > 0$ such that 
    $\Sigma^{*}$ is contained in $]\lambda_{1}^{m}+ \epsilon, + 
    \infty[ \times ]\lambda_{1}^{m} + \epsilon, + \infty[$ if and 
    only if $m$ has compact support in $]T_{1}, T_{2}[$.
\end{coro}

\begin{rem}\label{rem:5.6}
    {\rm When there is no weight, it was observed in \cite{DF-Go}, 
    \cite{A-C-G}, under various boundary conditions, that a strong 
    connexion (qualitative and quantitative) exists between on one 
    side the asymptotic behaviour of the first curves of $\Sigma^{*}$ 
    and on the other side the uniformity or nonuniformity of the 
    antimaximum principle~: the existence of an $\epsilon > 0$ as in 
    corollary \ref{coro:5.5} corresponds to the uniformity of the 
    antimaximum principle, and moreover the largest $\epsilon$ admissible corresponds exactly to the largest interval of 
    uniformity. In this context corollary \ref{coro:5.5} should be 
    compared with the recent result of \cite{G-G-P} where it is shown 
    that in the Dirichlet problem, whatever the weight, the 
    antimaximum principle is always non uniform. It follows that the 
    connexion referred to above does not hold anymore in the presence 
    of a weight with compact support, even $\geqs{\not \equiv} 0$.
    }
\end{rem}


\begin{thebibliography}{xx}

\bibitem{Al1} M. ALIF, Spectre de Fu\v cik : probl\`eme avec poids en dimension un et quelques remarques
en dimension sup\'erieure, Th\`ese, Universit\'e Libre de Bruxelles, 1999.

\bibitem{Al2} M. ALIF, Fu\v cik spectrum for the Neumann problem with indefinite weights, to appear.

\bibitem{star} M. ARIAS, J. CAMPOS, M. CUESTA and J.-P. GOSSEZ, Asymmetric elliptic problems with weights, to appear.

\bibitem{A-C-G}  M. ARIAS, J. CAMPOS and J.-P. GOSSEZ, On the
	antimaximum principle and the Fu\v cik spectrum for the Neumann
	$p$-laplacian, Diff. Int. Equat., 13 (2000), 217-226.


\bibitem{Ca} J. CAMPOS, Espectro de Fu\v cik para operadores elipticos, Th\`ese, Universidad de Granada, 1996.

\bibitem{Da} N. DANCER, On the Dirichlet problem for weakly nonlinear elliptic partial differential equations, Proc.
Royal Soc., Edimb., 76(1977), 283-300.

\bibitem{DF-Go} D. DE FIGUEIREDO and J.-P. GOSSEZ, On the first 
	curve of the Fu\v cik spectrum of an elliptic operator, Diff. Int. 
	Equat., 7 (1994), 1285-1302.  

\bibitem{Dr}  P. DRABEK, Solvability and bifurcations of nonlinear 
	equations, Pitman Research Notes in Mathematics, 264 (1992).

\bibitem{Fu} S. FU\v CIK, Solvability of nonlinear equations and boundary value problems, Reidel, Dordrecht, 1980.

\bibitem{G-G-P} T. GODOY, J.-P. GOSSEZ and S. PASZKA, Antimaximum principle for elliptic problems with weight, Electr. J.
Diff. Equat., 1999 (1999), 1-15.

\bibitem{Ha} P. HARTMAN, Ordinary differential equation, Wiley, 1964.

\bibitem{Ry} B. RYNNE, The Fu\v cik spectrum of general Sturm-Liouville problems, J. Diff. Equat., to appear.
\end{thebibliography}
\end{document}